\documentclass[12pt]{amsart}
\usepackage{amsmath,amssymb}
\usepackage{supertabular}
\usepackage{tikz}
\usepackage{multirow}
\usepackage{enumerate}

\usetikzlibrary{calc, shapes, backgrounds }
\usepackage{pgfplots,pgfplotstable}
\pgfplotsset{width=7cm}

\newtheorem{theorem}{Theorem}
\newtheorem{lemma}{Lemma}

\newtheorem{corollary}{Corollary}

\newtheorem{remark}{Remark}

\newtheorem{notation}{Notation}

\DeclareMathOperator{\mex}{mex}

\newcommand{\raf}[1]{(\ref{#1})}

\newcommand{\ga}{\alpha}
\newcommand{\gb}{\beta}

\def\G{\mathcal{G}}
\def\cG{\mathcal{G}}
\def\P{\mathcal{P}}

\def\Z{\mathbb{Z}}
\def\ZZ{\mathbb{Z}}

\def\eps{\epsilon}
\def\ge{\eta}

\def\cH{\mathcal{H}}

\begin{document}

\title{On the Sprague-Grundy function of {\sc Exact $k$-Nim}}

\author{Endre Boros}
\address{MSIS and RUTCOR, RBS, Rutgers University,
100 Rockafeller Road, Piscataway, NJ 08854}
\email{endre.boros@rutgers.edu}

\author{Vladimir Gurvich}
\address{MSIS and RUTCOR, RBS, Rutgers University,
100 Rockafeller Road, Piscataway, NJ 08854; \\
Dep. of Computer Sciences, National Research University,
Higher School of Economics (HSE), Moscow}
\email{vladimir.gurvich@rutgers.edu}

\author{Nhan Bao Ho}
\address{Department of Mathematics and Statistics, La Trobe University, Melbourne, Australia 3086}
\email{nhan.ho@latrobe.edu.au, nhanbaoho@gmail.com}

\author{Kazuhisa Makino}
\address{Research Institute for Mathematical Sciences (RIMS)
Kyoto University, Kyoto 606-8502, Japan}
\email{e-mail:makino@kurims.kyoto-u.ac.jp}

\author{Peter Mursic}
\address{MSIS and RUTCOR, RBS, Rutgers University,
100 Rockafeller Road, Piscataway, NJ 08854}
\email{peter.mursic@rutgers.edu}



\subjclass[2000]{91A46}
\keywords{Moore's {\sc Nim}, Exact {\sc Nim}, impartial combinatorial games, $\P$-position, Sprague-Grundy function.}

\thanks{Part of this research was done at the Mathematisches Forschungsinstitut Oberwolfach
during a stay within the {\it Research in Pairs} Program in 2015.
This research was partially supported by the Scientific
Grant-in-Aid from Ministry of Education, Science, Sports and Culture of Japan.
The second author was also supported by
the Russian Academic Excellence Project '5-100'.}

\begin{abstract}
Moore's generalization of the game of {\sc Nim} is played as follows.
Let $n$ and $k$ be two integers such that  $1 \leq k \leq n$.
Given $n$ piles of tokens, two players move alternately,
removing tokens from at least one and at most  $k$  of the piles.
The player who makes the last move wins.
The game was solved by Moore in 1910 and
an explicit formula for its Sprague-Grundy function
was given by Jenkyns and Mayberry in 1980, for the case $n = k+1$ only.
We introduce another generalization of {\sc Nim}, called
{\sc Exact $k$-Nim}, in which each move reduces exactly $k$ piles.
We give an explicit formula for the Sprague-Grundy function of
{\sc Exact $k$-Nim} in case  $2k \geq n$.
In case $n=2k$ our formula is surprisingly similar to Jenkyns and Mayberry's one.
\end{abstract}

\maketitle

\section{Introduction} \label{Ss.Intro}

We consider {\it combinatorial games} of two players;
they take turns alternating and one who makes the last move wins.
Both players have perfect information and there are no moves of chance.
A game is called {\it impartial} if both players
have the same possible moves in each position and
\emph{acyclic} if it is impossible to revisit the same position.
In this paper we consider only impartial acyclic combinatorial games
and call them simply games.
A more detailed introduction to combinatorial games can be found in \cite{BCG01-04, Con76}.

If there is a move from position $x$ to $y$, we write $x \to y$.
For a set $S$ of nonnegative integers the {\it minimum excluded value} of $S$
is defined as the smallest nonnegative integer that is not in $S$ and
is denoted by  $\mex (S)$. In particular, $\mex(\emptyset) = 0$.
The {\it Sprague-Grundy} (SG) value of a position $x$
in a game $G$ is defined recursively by  $$\G(x) = \mex \{g(y) \mid x \to y\}.$$
A position of SG value $t$ is called a {\it $t$-position};
$0$-positions are also known as $\P$-positions.
A player who moves into a $\P$-position can win the game.
The SG function is instrumental in the theory of disjunctive sums of games;
see \cite{BCG01-04, Con76, Gru39, Spr35, Spr37}.


A classical example is {\sc Nim} studied by Bouton \cite{Bou901}.
A position in {\sc Nim} consists of $n$ piles of tokens.
Two players alternately choose one of the piles and remove
an arbitrary (positive) number of tokens from that pile.
Bouton characterized the $\P$-positions and in fact described the SG function of {\sc Nim}.

Moore \cite{Moo910} introduced a generalization in which
a player can remove tokens from at least one and at most  $k$  of the piles, for some fixed  $k<n$.
We call this game Moore's {\sc Nim} and denote it by {\sc Nim}$^\leq_{n, k}$.
When $k=1$ {\sc Nim}$^\leq_{n, 1}$ is the traditional {\sc Nim}.

Moore generalized Bouton's results and
characterized the $\P$-positions of {\sc Nim}$^\leq_{n, k}$ as follows.
We denote by $\Z_\geq$ the set of nonnegative integers and
by $x \in \Z^n_\geq$ a position of the game,
where $x_i$ is the number of tokens in pile $i$.
Let us represent the components of $x$ in binary form
$x_i=\sum_{j}x_{ij}2^j$, $i=1,...,n$, define $y_j=\sum_{i=1}^n x_{ij} \mod (k+1)$, and set
\[
M(x)=\sum_{j} y_j (k+1)^j.
\]
For instance, if  $n=3$, $k=2$, and $x=(2,3,6)$, then using $3$ binary digits we can write:
$2=010$, $3=011$, and $6=110$ yielding
$y_0=1\mod 3=1$, $y_1=3\mod 3=0$, and $y_2=1\mod 3=1$,
from which we get $M(x)=1\cdot 3^0+0\cdot 3^1+1\cdot 3^2=10$.

Moore proved that $x$ is a $\P$-position of {\sc Nim}$^\leq_{n, k}$ if and only if  $M(x)=0$.
Berge \cite[Theorem 3, page~55]{Ber62} claimed that  $\G(x)$ is simply equal to  $M(x)$.
However, Jenkyns and Mayberry \cite{JM80} pointed out that this is an overstatement and
the equality holds only when  $M(x) \leq 1$ or $\G(x)\leq 1$.
For example, direct calculations show that
$2 = \G(x) < M(x) = 3$   for  $k = 2, n = 3, x = (0,0,2)$ and
$8 = \G(x) > M(x) = 2$   for  $k = 2, n = 3, x = (2,3,3)$.

For the case of $n = k+1$ Jenkyns and Mayberry \cite{JM80}
provided a formula for the SG function of {\sc Nim}$^\leq_{k+1, k}$.
An alternative proof for a slightly more general game
was given recently in \cite{BGHM15}.

In this paper we introduce another generalization of {\sc Nim}.
Given positive integers $n$  and  $k$  such that $1 \leq k \leq n$,
we define {\sc Exact $k$-Nim}, denoted by {\sc Nim$^=_{n,k}$},  as follows.
Given $n$ piles of tokens, by one move a player chooses exactly $k$ piles and
removes arbitrary positive number of tokens from each of them.
The game terminates when there are less than  $k$  nonempty piles.
{\sc Nim}$^=_{n, k}$  turns into the standard {\sc Nim} when  $k=1$ and
it is the trivial one-pile {\sc Nim} when $k=n$.

\subsection*{Main results}

Given a position $x\in \Z_\geq^n$ of {\sc Nim$^=_{n,k}$},
we denote by $T_{n,k}(x)$ the maximum number of consecutive moves
one can make starting with  $x$.
We call $T_{n,k}$ the \emph{Tetris} function of the game.

The following two theorems characterize the SG function of {\sc Nim$^=_{n,k}$} for $2k\geq n$.
\begin{theorem} \label{thm.n<2k}
If $n < 2k$ then the SG function of {\sc Nim$^=_{n,k}$} is
equal to its Tetris function,  $\G(x) = T_{n,k}(x)$.
\end{theorem}

\begin{theorem} \label{thm.n=2k}
Let $k \geq 2$, $n = 2k$, and
let $x = (x_1, \ldots, x_n)$ be a position of {\sc Nim}$^=_{2k, k}$.
Set
\begin{align}
u(x) &= T_{2k,k}(x),  \label{eq:u}\\
m(x) &= \min_{1\leq i\leq 2k} x_i,\label{eq:m}\\
y(x) &= T_{2k,k}\big(x_1 - m(x), \ldots, x_{2k} - m(x)\big),\label{eq:y}\\
z(x) &= 1+\binom{y(x)+1}{2}, \text{ and}\label{eq:z}\\
v(x) &= \big(z(x)-1\big) +\big[\big(m(x)-z(x)\big)\mod \big(y(x)+1\big)\big].\label{eq:v}
\end{align}
Then the SG function of {\sc Nim}$^=_{2k,k}$ is given by formula
\begin{equation}\label{eq-binomial;n=2k-Formula}
\G(x) ~=~ \begin{cases}
u(x), & \text{if~~}  m(x) < z(x);\\
v(x), & \text{if~~}  m(x) \geq z(x).
\end{cases}
\end{equation}
\end{theorem}

Note that this formula fails for  $k = 1$.
In this case {\sc Nim}$^=_{2, 1}$ is the standard $2$-pile {\sc Nim} and
its SG function is the modulo $2$ sum of the two coordinates of a position,
as described by Bouton. 
This function is different from the one described by the above formula.

We also would like to remark that the above formula is surprisingly
similar to the one given by Jenkyns and Mayberry
in \cite{JM80} for {\sc Nim}$^\leq_{k+1, k}$.

\bigskip

The above result implies a simple characterization of
the $0$- and $1$-positions of {\sc Nim}$^=_{2k, k}$.
A position $x = (x_1, \ldots, x_n)$ is said to be
\emph{nondecreasing} if $x_1 \leq x_2 \leq \cdots \leq x_n$.

\begin{corollary}   \label{cor-0,1-exact:n=2k}
Given a nondecreasing position $x = (x_1, \dots, x_{2k})$ of the game {\sc Nim}$^=_{2k, k}$,
\begin{enumerate} \itemsep0em
\item [\rm{(i)}] $x$  is a $\P$-position if and only if
more than half of its smallest coordinates are equal, that is,  $x_1 = \cdots = x_k = x_{k+1}$.
\item [\rm{(ii)}] $x$  is a $1$-position if and only if $x_1 = \cdots = x_{k - \ell} = 2c$  and
$x_{k-\ell+1} = \cdots = x_{k+ \ell + 1} = 2c+1$
for some integer  $c \in \ZZ_\geq$  and  $\ell \in \{0,1, \ldots, k-1\}$.
\end{enumerate}
\end{corollary}

Both statements \rm{(i)} and  \rm{(ii)} follow from Theorem \ref{thm.n=2k},
but can also be derived much simpler, directly from the definitions.

\medskip

The case of $2k < n$ looks much more difficult and it is still open.
Moreover, we have not even been able to characterize
the $\P$-positions of {\sc Exact $k$-Nim} for $1 < k < n/2$,
e.g., for {\sc Nim}$^=_{5,2}$.

\bigskip

The rest of the paper is organized as follows.
In Section \ref{Ss.n<2k} we characterize the SG function of {\sc Nim$^=_{n,k}$} for the case $2k>n$.
In Section \ref{Ss.n=2k} we characterize the SG function of {\sc Nim$^=_{n,k}$} for the case $2k=n$.
In Section \ref{Ss.Moore01} we provide an alternative proof for the above stated result of Jenkyns and Mayberry \cite{JM80}.
Finally in Section \ref{Ss.Tetris} we show that for a given position we can compute efficiently the corresponding SG value.


\section{SG function in the case of $n < 2k$} \label{Ss.n<2k}

For our proof we need the following basic properties of the Tetris function.
Given positions $x,x'\in\Z^n_\geq$ we write $x\leq x'$ if $x_i\leq x'_i$ holds for $i=1,...,n$.

\begin{lemma}\label{tbasic}
Consider two positions $x, x'\in\ZZ^n_\geq$.
\begin{itemize}
\item [\rm{(i)}] If $x'\leq x$ then $T_{n,k}(x')\leq T_{n,k}(x)$.
\item [\rm{(ii)}] If in addition we have $\sum_{i=1}^n(x_i-x_i')=1$, then $T_{n,k}(x)-1\leq T_{n,k}(x')\leq T_{n,k}(x)$.
\end{itemize}
\end{lemma}
\proof
It is immediate by the definition.
\qed

A move in {\sc Nim}$^=_{n, k}$ is called {\it slow} if exactly one token is taken from each of the $k$ chosen piles.



\begin{lemma}\label{Tetrisnotdecrease}
Consider a position $x=(x_1, \ldots, x_n) \in \Z^n_\geq$ with some indices $i, j$ such that
$x_i < x_j$. Let $x'=(x'_1, \ldots, x'_n)$ be defined by
\begin{equation}\label{eq-l2}
x_l'=
\begin{cases}
x_i+1, & \text{ if } l=i, \\
x_j-1, & \text{ if } l=j, \\
x_l,   & \text{ otherwise. }\\
\end{cases}
\end{equation}
Then we have $T_{n,k}(x) \leq T_{n,k}(x')$.
\end{lemma}
In other words, the Tetris function is nondecreasing when
we move a token from a larger pile to a smaller one.

\proof
Consider any sequence of slow moves from $x \to \cdots \to x''$. If $x''_j>0$  then the same sequence of slow moves can be made from $x'$ since $x'_l \geq x_l$ for $l \neq j$.

If $x''_j=0$ then since $x_j > x_i$, this sequence contains a slow move reducing $x_j$ but not $x_i$.
Let us modify this move reducing $x_i$ rather than $x_j$ and keeping all other moves of the sequence unchanged.
The obtained sequence has the same length and consists of slow moves from $x'$. \qed

Notice that we can generalize  Lemma \ref{Tetrisnotdecrease}
replacing $\pm 1$ in \raf{eq-l2} by $\pm \Delta$ for any integral $\Delta \in [0,x_j-x_i]$.

\begin{lemma}\label{bestslowmove}
The slow move that reduces the $k$ largest piles of $x$ reduces the Tetris value $T_{n,k}(x)$ by exactly one.
\end{lemma}

\proof
Let $x'$ be the position obtained from $x$ by reducing the $k$ largest piles of $x$ by exactly one each. Let $x''$ be another position obtained by some slow move.
By applying \raf{eq-l2} repeatedly, we can obtain $x'$ from $x''$ with $T_{n,k}(x'') \leq T_{n,k}(x')$ by Lemma \ref{Tetrisnotdecrease}.
This implies that $x'$ has the  highest Tetris value among all positions each reachable from x by a slow move. By Lemma \ref{tbasic},
each slow move reduces the Tetris value by at least one and there exists a slow move reducing it by exactly one.
Hence, $T_{n,k}(x')=T_{n,k}(x)-1$. \qed

\bigskip
\noindent\textbf{Proof of Theorem }\ref{thm.n<2k}:
The Tetris value is the largest number of moves one can take from a position $x$,
implying that $\G(x)$ is at most $T_{n,k}(x)$.
Therefore, it is enough to show that for all integral $g$ such that
$0 \leq g < T_{n,k}(x)$ there exists a move $x\to x'$ such that $T_{n,k}(x')= g$.

Consider the move $x \to x'$ that reduces the largest $k$ piles to $0$.
For the resulting position we have $T_{n,k}(x')=0$ because $2k>n$.
Let us also consider the move $x \to x''$ that reduces the $k$ largest piles each by only $1$.
Then we have $T_{n,k}(x'')=T_{n,k}(x)-1$ according to Lemma \ref{bestslowmove}.
Any position between $x'$ and $x''$ is reachable form $x$. Thus, by Lemma \ref{tbasic} the claim follows.  \qed

\section{SG function in case of $n=2k$} \label{Ss.n=2k}

Let us consider the game {\sc Nim}$^=_{2k, k}$, where $k\geq 2$, and let $x$ be a position of this game.
Recall that to $x$ we associated several parameters in \eqref{eq:u}--\eqref{eq:v}.
Based on these parameters, we classify the positions into the following two types:
\begin{itemize} \itemsep0em
\item [\rm{(i)}]   \emph{type I}, if $m(x)<z(x)$, and
\item [\rm{(ii)}]  \emph{type II} otherwise.
\end{itemize}

We shall need some technical lemmas for our proof. In this section, all positions belong to the game {\sc Nim}$^=_{2k, k}$.

\begin{lemma} \label{continuity}
Consider a position $x$ and two moves $x\to x'$ and $x\to x''$ such that $x'\geq x''$ $($componentwise$)$ and $y(x')\geq y(x'')$.
Then, for every integer $g$ with $y(x')\geq g \geq y(x'')$ there exists a move $x\to x'''$ such that $y(x''')= g$ and $x'\geq x'''\geq x''$.
\end{lemma}
\proof
Note that $x'\geq x''$ implies that in the two moves $x\to x'$ and $x\to x''$ the same $k$ components are decreased.
Let us denote by $K\subseteq \{1,2,...,2k\}$ these components.
Let us now start decreasing the components $x'_j$, $j\in K$ one by one,
keeping their values always greater than the corresponding $x''_j$ values.
After $\sum_{j\in K} x'_j-x''_j$ steps we reach $x''$.
In each of these steps the corresponding value $y$ can only decrease, and by at most $1$.
Hence there will be at least one such position $x'\geq x'''\geq x''$ with $y(x''') = g$.
It also follows that $x\to x'''$ is a move, completing the proof.
\qed

Let us also note that in fact $z$ defined in \eqref{eq:z} depends uniquely on $y(x)$ and hence can be considered as a function of $y$.
We shall also need the following easy arithmetical facts.
\begin{lemma}\label{interval}
Every nonnegative integer $g$ belongs to exactly one of the intervals
$$[z(y)-1, z(y)+y-1] = \left[\binom{y+1}{2},\binom{y+2}{2}-1\right] \textrm{ for }y\in \{0, 1, 2, ...\}.$$
\qed
\end{lemma}
\begin{corollary}\label{delta-epsilon}
For every nonnegative integer $g$, there exist unique $\nu(g)$ and $\eps(g)$ integer values such that
\[
g ~=~ \binom{\nu(g)+1}{2} +\eps(g) ~~~\textbf{ and }~~~ 0\leq \eps(g)\leq \nu(g).
\]
\end{corollary}

\smallskip

Given a nondecreasing position $x$, let us construct another position $\bar x$ from $x$ by emptying the first $n-k$ piles and adding these $\sum_{i=1}^{n-k}x_i$ tokens, one by one, to the last $k$ piles as follows:
In each step we add one token to the smallest of these $k$ piles.
If there are several such piles, we break the tie by adding this token to the pile of the largest index.

It is easy to see that we have $T_{n,k}(\bar x)=\min(\bar x_i  |  n-k+1 \leq i \leq n)=\bar x_{n-k+1}$.

\begin{figure}[ht]
\centering
\includegraphics[height=6cm]{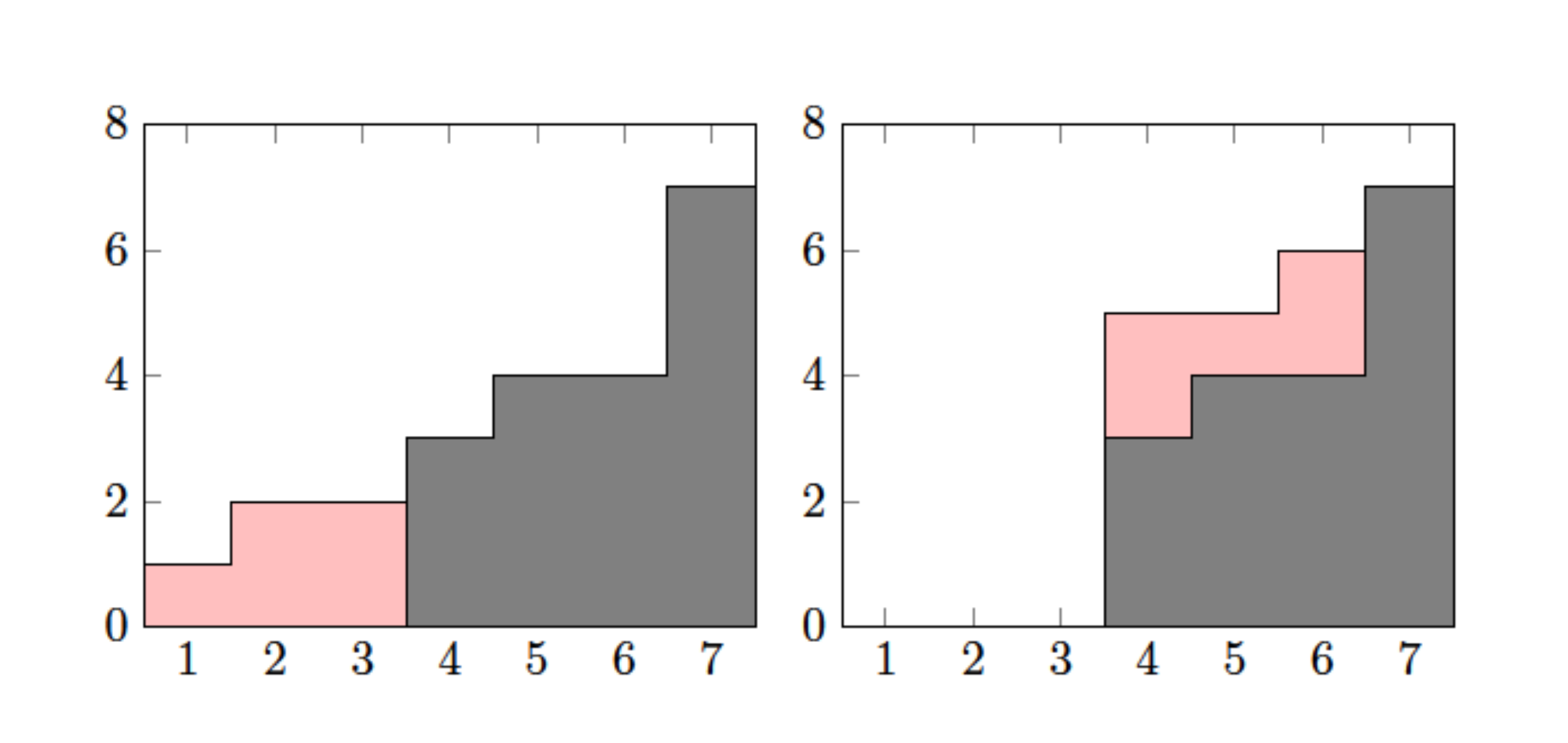}
\caption{$k=4,n=7$, $x=(1,2,2,3,4,4,7)$, and $\bar x=(0,0,0,5,5,6,7)$.} \label{ADD(1)}
\end{figure}

\begin{lemma}\label{txequalstxbar}
The above construction of $\bar x$ from $x$ keeps the Tetris value unchanged: $T_{n,k}(\bar x)=T_{n,k}(x)$.
\end{lemma}
\proof
Let us note that $T_{n,k}(\bar x) \leq T_{n,k}(x)$ by Lemma \ref{Tetrisnotdecrease}.
By the above definition of $\bar x$, none of the tokens from the smallest $n-k$ piles of $x$ are moved to any pile of size larger than $T_{n,k}(\bar x)+1$ and hence we have
\begin{equation}\label{volumeargument}
V(x):=\sum_{i=1}^n \min(x_i,T_{n,k}(\bar x)+1) =\sum_{i=n-k+1}^n \min(\bar x_i,T_{n,k}(\bar x)+1).
\end{equation}
Since $\bar x_{n-k+1}=T_{n,k}(\bar x)$, we get by \eqref{volumeargument} that
\begin{equation*}
V(x)  \leq k-1+ k T_{n,k}(\bar x) < k(T_{n,k}(\bar x)+1).
\end{equation*}
Assume now indirectly that $T_{n,k}(\bar x) < T_{n,k}(x)$.
Then it is possible to construct a  sequence of $T_{n,k}(\bar x)+1$ slow moves from $x$.
By such sequence any pile would be reduced at most $T_{n,k}(\bar x)+1$ times,
and therefore the total number of the removed tokens is at most $V(x)$,
implying $k(T_{n,k}(\bar x)+1)\leq V(x)$, contradicting the above inequality.
The obtained contradiction implies that $T_{n,k}(x)=T_{n,k}(\bar x)$. \qed

\begin{lemma}\label{getmeTetris}
Let $x = (x_1, \ldots, x_n)$ be a position.
If there exists a nonnegative integer $g$ such that
\begin{equation}\label{eq-getmeTetris}
k g \leq \sum_{i=1}^{n} \min(x_i,g) \quad \text{ and } \quad \sum_{i=1}^{n} \min(x_i,g+1) < k(g+1),
\end{equation} \label{E.delta}
then we have $T_{n,k}(x)=g$.
\end{lemma}
\proof
The sequence $s(g)=\frac{1}{g}\sum_{i=1}^n \min(x_i,g)$ is monotone non-increasing for $g \in \ZZ_\geq$,
and hence, the inequalities \raf{eq-getmeTetris} can hold for at most one $g$.
Without any loss of generality we can assume that $x$ is nondecreasing, and thus we can define $\bar x$, as above.
Lemma \ref{txequalstxbar} then implies that \raf{eq-getmeTetris} holds for $g=T_{n,k}(x)=T_{n,k}(\bar x)$. \qed

\bigskip

\noindent{\bf Proof of Theorem \ref{thm.n=2k}.}
Our main result claims that the SG function $\G(x)$ of
{\sc Nim}$^=_{2k, k}$ for $k\geq 2$ is equal to the function

$$
g(x)=\begin{cases}
u(x), & \text{if }m(x) < z(x); \\
v(x), & \text{if }m(x) \geq z(x);
\end{cases}
$$
where $ u,m,y,v,$ and $z$ are defined in \eqref{eq:u}--\eqref{eq:v}.
To prove this theorem, it is enough to show the following two properties of $g(x)$:
\begin{itemize}
\item[(I)] for any move $x \rightarrow x'$ we have $g(x) \neq g(x')$ and
\item[(II)] for every value $\delta$ such that $0 \leq \delta < g(x)$ there exists a move $x \to x'$ with $g(x')=\delta$.
\end{itemize}

We assume that $x$ is nondecreasing.

\subsection{Proof of (I)}

We now prove property (I). Let us first consider the case $z(x) > m(x)$.
Obviously, any move $x \rightarrow x'$ reduces the Tetris value by at least 1, implying $u(x) > u(x')$.
Using this and the definitions, we get $g(x)=u(x) > u(x')$. If $g(x') = u(x')$ then by the above inequality, we get $g(x) \neq g(x')$.
On the other hand, if $g(x') = v(x')$ then we  must have $m(x') \geq z(x')$ and thus $v(x')=z(x')-1 + \left((m(x')-z(x')) \mod (y(x')+1)\right)\leq z(x')-1 + (m(x')-z(x'))=m(x')-1$.
Thus, we get $g(x) = u(x) > u(x') > m(x')-1 \geq  v(x') = g(x')$.

It remains to consider the case $z(x) \leq m(x)$, in which case $v(x)\leq m(x)-1$ follows by the definitions.  
	\begin{enumerate} \itemsep0em
	\item Suppose $z(x') > m(x')$. 
    We can estimate $u(x') \geq m(x)+m(x')$ since $2k = n$.
    Furthermore, $g(x)=v(x) \leq m(x)-1$ and thus $g(x')=u(x') \geq m(x)+m(x') > m(x)-1 \geq g(x)$.
    \item Suppose $z(x') \leq m(x')$. Then $g(x')=v(x')$.  Note that $z(x) \leq m(x)$ and thus $g(x)=v(x)$.
    Also note that $m(x') \leq m(x)$. We examine the last inequality.
    	\begin{enumerate} \itemsep0em
	    \item Suppose $m(x')=m(x)$.
                Then $x-m(x)\to x'-m(x')$ is a move, and hence decreases the Tetris value implying $y(x)>y(x')$.
	           By Lemma \ref{interval}, we have
                $z(x')+y(x')-1 < z(x)-1$
                and so
                $z(x)  > z(x') + y(x')$
                implying 
                $v(x) \geq z(x)-1> z(x')-1+y(x') \geq v(x')$,
                since
                $y(x') \geq \bigg(\big(m(x') - z(x)\big) \mod\big(y(x) + 1\big)\bigg)$,
                regardless  of the value of $m(x') - z(x)$.
	    \item Suppose $m(x') < m(x)$. We compare $y(x)$ with $y(x')$.
            \begin{enumerate} \itemsep0em
	        \item If $y(x)= y(x')$ then $z(x)=z(x')$. By the definition of a legal move, $x'$ has at least $k$ piles not smaller than $m(x)$. Therefore, $1 \leq
            m(x)-m(x') \leq y(x')=y(x)$, which implies that $m(x) \mod (y(x)+1) \neq m(x') \mod (y(x)+1)$ and thus $v(x) \neq v(x')$.
	        \item If $y(x) \neq y(x')$, since the $y$ values are different, by Lemma \ref{interval} we have that $v(x)$ and $v(x')$ are in different intervals,
            therefore $v(x)\neq v(x')$.
            \end{enumerate}
	   \end{enumerate}
	\end{enumerate}

\subsection{Proof of (II)}

We prove property (II) by considering type I and type II positions, separately.

\subsubsection{Type I positions: $m(x)<z(x)$}

First, let us consider the case $m(x)=0$. Then there are at most $2k-1$ nonempty piles.
So we can reduce the Tetris value $T_{2k,k}(x)$ to $0$ by emptying the $k$ largest piles of $x$.
Therefore, by Lemmas \ref{achieveTetrisvalue} and \ref{continuity}, for any $0 \leq \delta < T_{2k,k}(x)$, there exists a move $x \to x'$ with $T_{2k,k}(x')=\delta$.
All such moves also have $m(x')=0$, therefore $m(x')<z(x')$.

From now on we can assume that $m(x) \geq 1$. Let us consider the following four subcases, depending on the value $\delta$.

\begin{enumerate}
\item $0\leq \delta <  m(x)$.
To simplify our proof, let us use simply $m$ instead of $m(x)$ in this section.
Let us observe first that since $x$ is type I we have $m<z(x)$ implying by \eqref{eq:z} that
\begin{equation}\label{sqrt-m1}
y(x)>\nu(m-1)\geq 0.
\end{equation}
with $\nu$ being defined as in Corollary \ref{delta-epsilon}.

Let us next define a set $Q=Q(x)$ of pairs of integers by setting
\[
Q^1(x)=\left\{(\mu,\ge)\left| \begin{array}{c}m-\ge~\leq~ \mu~\leq~ m\\
\ge ~\leq \nu(m-1)-1
\end{array}\right.\right\}
\]
\[
Q^2(x)=\left\{(\mu,\ge) \left| \begin{array}{c}m-\eps(m-1)~\leq~ \mu~\leq~ m\\
\ge ~= \nu(m-1)
\end{array}\right. \right\} ,
\]
and defining $Q=Q^1(x)\cup Q^2(x)$.

We show that if for a position $x^*$ we have $(m(x^*),y(x^*))\in Q$, then $x^*$ is of type II.
To see this consider a pair $(\mu,\ge)\in Q^1(x)$.
Then we have by the definition of $Q^1(x)$ that $\mu\geq m-\ge$ and that $\nu(m-1)-1\geq \ge$ from which $z(\nu(m-1)-1)\geq z(\ge)$ follows by the definition of $z$ in \eqref{eq:z}.
We also have the inequality $m-\nu(m-1)+1\geq z(\nu(m-1)-1)$ since $m\geq \binom{\nu(m-1)+1}{2}$ by the definition of $\nu$ in Corollary \ref{delta-epsilon}.
Putting these together, we obtain $\mu\geq z(\ge)$ as stated.
For $(\mu,\ge)\in Q^2(x)$ we have $\mu\geq m-\epsilon(m-1)$ and $\ge=\nu(m-1)$ by the definition of $Q^2(x)$.
Since $m-1=\binom{\nu(m-1)+1}{2}+\epsilon(m-1)$ by Corollary \ref{delta-epsilon},
the inequality $\mu\geq m-\epsilon(m-1)=z(\ge)$ follows again.

We show next that for all pairs $(\mu,\ge)\in Q$ there exists a move $x\to x^*$ such that $\mu=m(x^*)$ and $\ge=y(x^*)$ (and $x^*$ is of type II, as we argued in the previous paragraph.)
For this let us consider first $\mu\leq m-1$ and note that if $(\mu,\ge)\in Q$ then $m-\mu\leq \ge\leq \nu(m-1)$.
Our plan is to use Lemma \ref{continuity} and to cover this range of $\ge$ values by two constructions.

Let us define a pair of positions $x'\geq x''$ by
\begin{align*}
x'_i &=
\begin{cases}
x_i,   &\text{ for } i=1 \text{ and } i\geq k+2,\\
\mu,   &\text{ for } i=2,                    \\
x_i-1, &\text{ for } i=3,...,k+1
\end{cases} \\
\intertext{and}
x''_i &=
\begin{cases}
x_i,  &\text{ for } i=1 \text{ and } i\geq k+2,\\
\mu,  &\text{ for } i=2,...,k+1.
\end{cases}
\end{align*}
Note that since $\mu<m=x_1\leq x_2\leq x_i$ for $i\geq 3$, both $x'$ and $x''$ are reachable from $x$. We claim next that $y(x'')=m-\mu$ and $y(x')\geq \nu(m-1)$.
The first claim follows easily, since in $x''$ we have exactly $k$ positions larger than $\mu$, and $x_1=m=\mu +(m-\mu)$.
For the second let us note that since $x_1-\mu=m-\mu\geq 1$ we have $y(x')\geq x_{k+1}-\mu$. We can now apply Lemma \ref{continuity} for $x'$ and $x''$,
and conclude that for all values $m-\mu \leq \ge \leq x_{k+1}-\mu\leq y(x')$ there exists a move $x\to x^*$ such that $m(x^*)=\mu$ and $y(x^*)=\ge$.
For larger values of $\ge$ we need a modified construction:
\begin{align*}
x'_i & =
\begin{cases}
\mu&\text{ for } i=1,\\
x_i-1, &\text{ for } i=2,...,k,\\
x_i,   &\text{ for } i=k+1,...,2k,
\end{cases}\\
\intertext{and}
x''_i & = \begin{cases}
\mu,   &\text{ for } i=1,...,k,\\
x_i,   &\text{ for } i=k+1,...,2k.
\end{cases}
\end{align*}
Assuming $\mu < m$, we have $y(x')\geq y(x)\geq \nu(m)\geq \nu(m-1)$,
while $y(x'')=x_{k+1}-\mu$. We have again $m(x')=m(x'')=\mu$. Thus by applying
Lemma \ref{continuity} for $x'$ and $x''$ we can conclude that for all values $x_{k+1}-\mu\leq \ge \leq \nu(m-1)$ there exists a move $x\to x^*$ such that $m(x^*)=\mu$ and $y(x^*)=\ge$.

\medskip

Finally, for $\mu=m$ we proceed analogously, but with a third construction. Note first that $(m,\ge)\in Q$ if and only if $0\leq \ge\leq \nu(m-1)$.
Let us now proceed with constructing two positions reachable form $x$:
\begin{align*}
x'_i &=
\begin{cases}
x_i,  &\text{ for } i=1,...,k,\\
x_i-1, &\text{ for } i=k+1,...,2k
\end{cases} \\
\intertext{and}
x''_i &=
\begin{cases}
x_i,   &\text{ for } i=1,...,k,\\
m,     &\text{ for } i=k+1,...,2k.
\end{cases}
\end{align*}
It is easy to see that both are reachable from $x$, and that $m(x')=m(x'')=m$, $y(x')\geq y(x)-1\geq \nu(m-1)$ by \eqref{sqrt-m1} and $y(x'')=0$.
Thus, the existence of an $x^*$ reachable from $x$ with $m(x^*)=m$, $y(x^*)=\ge$ follows by Lemma \ref{continuity} for all $0\leq \ge\leq \nu(m-1)$.

By the above arguments we have a move $x\to x^*$ to a type II position $x^*$ with $m(x^*)=\mu$ and $y(x^*)=\ge$ for all $(\mu,\ge)\in Q$.
To conclude the proof for this case we claim that the corresponding $v(x^*)$ values include all integers in the interval $[0,m-1]$.
Note that $v(x^*)$ depends only on $m(x^*)$ and $y(x^*)$ for type II positions by \eqref{eq:v},
and that for a fixed value of $\ge\leq \nu(m-1)-1$ we have exactly $\ge+1$ consecutive integer values for $\mu$ such that $(\mu,\ge)\in Q$, implying that the corresponding values
$v(\mu,\ge)=\big(z(\ge)-1\big) + \big[\big(\mu-z(\ge)\big)\mod (\ge+1)\big]$
is exactly the set of integers in the interval
$[z(\ge)-1, z(\ge)+ \ge-1]$.
Thus our claim follows by the construction of $Q$ and by Lemma \ref{interval}.

\bigskip

\item $m(x) \leq \delta < x_2$ (Figure \ref{F1}).
Set $x'_i=0$ for $i=k+1, \ldots ,n-1$ and $x'_n=\delta-x_1$.

By definition $x_n\geq x_2>\delta\geq x_1$ which implies $x'_n=\delta-x_1<x_n-x_1\leq x_n$.
Furthermore, we have $x_i\geq x_2>0$ for $i=k+1,...,n-1$,
and therefore we indeed decrease exactly $k$ piles of $x$ to obtain $x'$. Thus $x'$ is reachable from $x$.

For $x'$ we have
\begin{align*}
\sum_{i=1}^{2k}\min(x_i',\delta)   &\geq x_1' + (k-1)\delta + (\delta-x_1') =k\delta \\
\intertext{and}
\sum_{i=1}^{2k}\min(x_i',\delta+1) &\leq x_1' + (k-1)(\delta+1) + (\delta - x_1)\\
                                   &= \delta + (k-1)(\delta +1) < k(\delta +1).
\end{align*}
Therefore, by Lemma \ref{getmeTetris}, we have $T_{2k,k}(x')=\delta$.

Since $k\geq 2$, $x'_{k+1}=0$ and therefore $m(x')=0$. Thus, $m(x')=0<1\leq z(x')$ implying
that $x'$ is of type I, from which $g(x')=u(x')=\delta$ follows by the above.
\begin{figure}[ht]
\centering
\includegraphics[height=7cm]{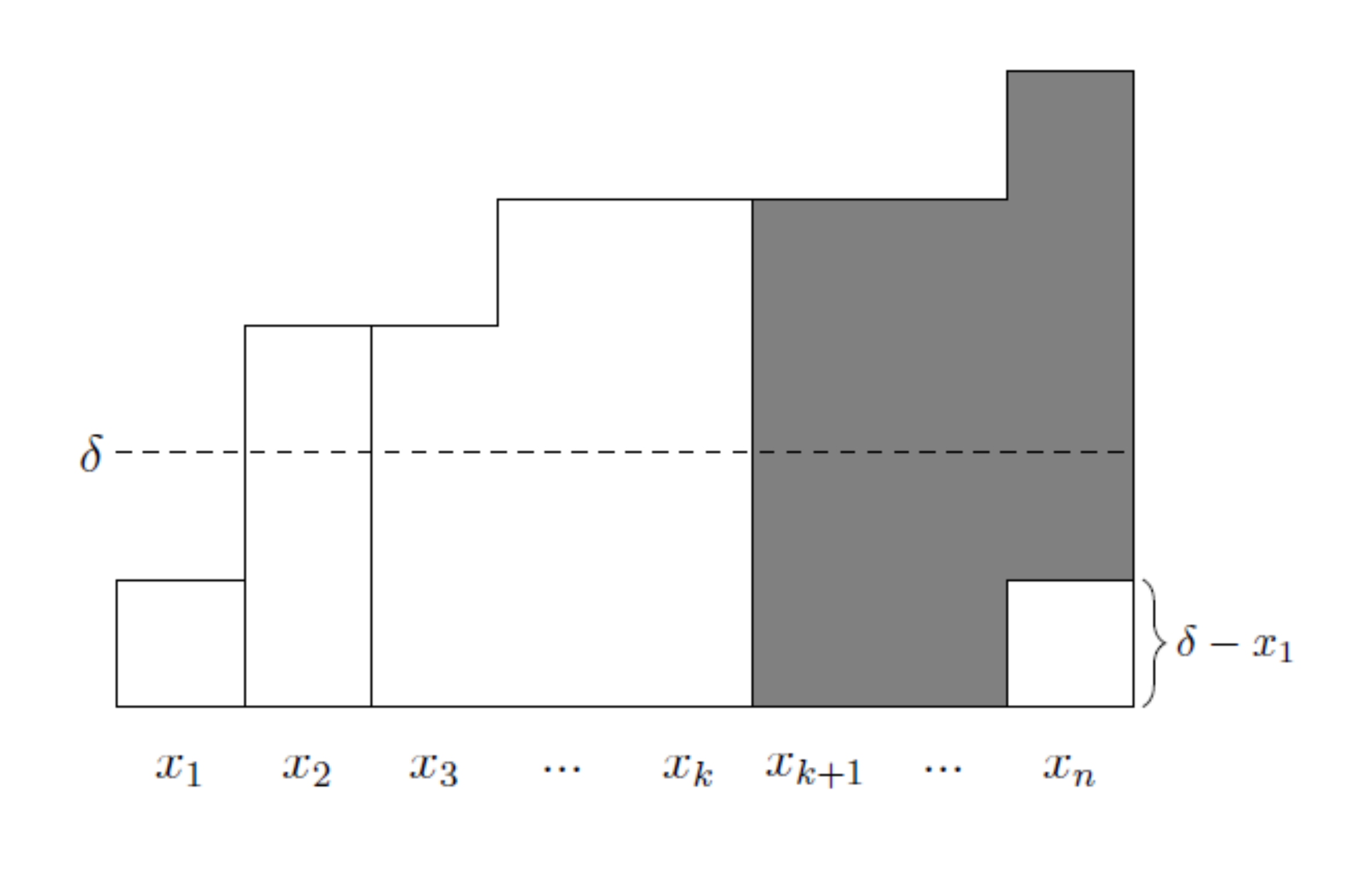}
\caption{$x'$ is obtained by removing the gray area.} \label{F1}
\end{figure}

\item $x_2 \leq \delta  < u(x)-m(x)$ (Figure \ref{F2}).  Let $A=\sum_{i=2}^{k+1} \max(0,\delta-x_i)$.
For $i=k+2, \ldots ,n$ choose $a_i$ such that $0 \leq a_i \leq \min(x_i-1,\delta)$ and $\sum_{i=k+2}^n a_i=A$.
First let us prove that this is possible, or equivalently that $A \leq \sum_{i=k+2}^n \min(\delta,x_i-1)$.
To see this let us consider two cases.
If $x_{k+1}>\delta$ then $\sum_{i=k+2}^n \min(\delta,x_i-1)\geq (k-1)\delta \geq A$.
If $x_{k+1}\leq \delta$ then let us define
\begin{equation} \label{B.def}
B=\sum_{i=2}^{k+1}x_i = k\delta - A,
\end{equation}
and observe that for any integer $t\geq\delta >x_{k+1}$ we have
\[
\sum_{i=1}^n\min(t,x_i) = m + B +\sum_{i=k+2}^n\min(t,x_i).
\]
Consequently, since we have $T_{2k,k}(x)=u(x)$, by Lemma \ref{getmeTetris} we can write for $t=u(x)$ that
\[
kt \leq \sum_{i=1}^n\min(t,x_i) = m(x)+B+\sum_{i=k+2}^n\min(t,x_i).
\]
Note that if we decrease $t=u(x)$ by $1$,
then the left hand side decreases by $k$ while the right hand side decreases by at most $k-1$,
hence the inequality remains valid. Let us repeat this $m(x)$ times, obtaining the inequality
\[
k(u(x)-m(x))+km(x) \leq m(x)+B+\sum_{i=k+2}^n\min(u(x)-m(x),x_i)+(k-1)m(x)
\]
from which
\[
k(u(x)-m(x))\leq B+\sum_{i=k+2}^n\min(u(x)-m(x),x_i)
\]
follows. Let us now decrease $t=u(x)-m(x)$ further by $1$, as well as replace $x_i$ by $x_i-1$.
Then the left hand side decreases by exactly $k$, while the right hand side decreases by at most $k$, yielding the valid inequality
\[
k(u(x)-m(x)-1)\leq B+\sum_{i=k+2}^n\min(u(x)-m(x)-1,x_i-1).
\]
Finally, we can decrease $t=u(x)-m(x)-1$ further on both sides to $t=\delta$ and similarly to the above argument obtain
\[
k\delta\leq B+\sum_{i=k+2}^n\min(\delta,x_i-1).
\]
By \eqref{B.def} we obtain the claimed inequality, and hence the proof for the existence of the $a_i$ values for $k=k+2,...,n$
that satisfy the desired inequalities.

Let us now consider the position $x'$ defined by
\[
x'_i =
\begin{cases}
0,     & \textrm{ for } i=1; \\
x_i,   &\textrm{ for }i=2, \ldots, k+1; \\
a_i,   &\textrm{ for }i=k+2, \ldots ,n.
\end{cases}
\]
By the above arguments $x\to x'$ is a move in the game.
The equality $g(x')=u(x')=\delta$ now follows by the above analysis and Lemma \ref{getmeTetris}, completing our proof in this case.

\begin{figure}[ht]
\centering
\includegraphics[height=9cm]{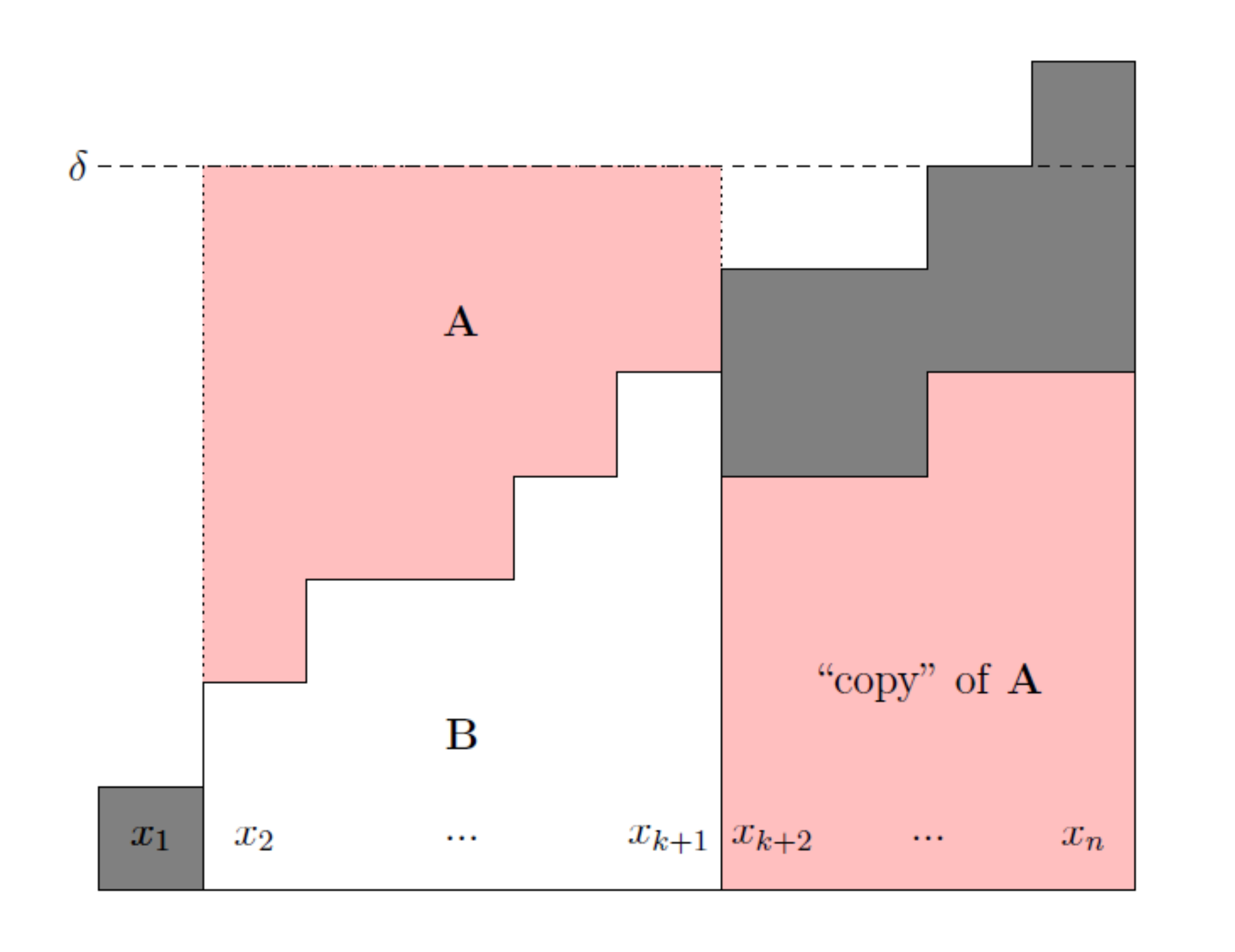}
\caption{$x'$ is obtained by removing the gray area.} \label{F2}
\end{figure}

\item $u(x)-m(x) \leq \delta < u(x)$ (Figure \ref{F3}).
Let us define position $x'$ as follow
$$
x_i' =
\begin{cases}
x_i-u(x)+\delta,  &\textrm{ for } i \in I_1=\{i \mid i\leq k, x_{i+k} \leq \delta\};   \\
x_i-u(x)+\delta,  &\textrm{ for } i \in I_2=\{i \mid i > k,   \delta < x_i \leq u(x)\}; \\
\delta ,           &\textrm{ for }i \in I_3=\{i \mid  i > k,   x_i > u(x) \}; \\
x_i,              &\textrm{ otherwise.}
\end{cases}
$$
Note that $x_i \geq m(x)\geq u(x)-\delta$, therefore $x_i'$ are all nonnegative.
It is easy to see that $I_1+k,I_2,I_3$ form a partition of $\{ k+1, \ldots , n \}$.
We have reduced $x_i$ for all $i \in I_1 \cup I_2 \cup I_3$, therefore $x\to x'$ is a move.

Next we note that the above construction implies that
\[
\sum_{i=1}^n \min(x_i',\delta) \geq \sum_{i=1}^n\min(x_i,u(x)) - k (u(x)-\delta)\geq k \delta
\]
where the last inequality follows by the fact that $u(x)=T_{2k,k}(x)$, and hence $\sum_{i=1}^n\min(x_i,u(x))\geq ku(x)$ by Lemma \ref{getmeTetris}.
Similarly, we get
\[
\sum_{i=1}^n \min(x_i',\delta+1) \leq \sum_{i=1}^n\min(x_i,u(x)+1) - k (u(x)-\delta)< k (\delta+1)
\]
since $\sum_{i=1}^n\min(x_i,u(x)+1)< k(u(x)+1)$.
Therefore $T_{n,k}(x')=u(x')= \delta$ follows by Lemma \ref{getmeTetris}.

Note that $x_{k+1}\leq \delta$, since otherwise $u(x)\geq \delta+m(x)$ would follow, contradicting our choice of $\delta$.
It follows that $x_1 \in I_1$ and thus $m(x')=m(x)-u(x)+\delta<m(x)$.

Let us define $\hat{x}_i=\min(x_i,u(x))$ for $i=1,...,n$. Then, by Lemma \ref{Tetrisnotdecrease}, we have $y(x)=y(\hat{x})$ and by construction of $x'$
we have $\hat{x}_i-m(x)\leq x'_i-m(x')$ for all indices $i$, implying $y(x)=y(\hat{x}) \leq y(x')$ by Lemma \ref{tbasic}.

Finally, $m(x')=m(x)-u(x)+\delta<m(x)<z(x)\leq z(x')$, which implies $g(x')=u(x')=\delta$.

\begin{figure}[ht]
\centering
\includegraphics[height=10.5cm]{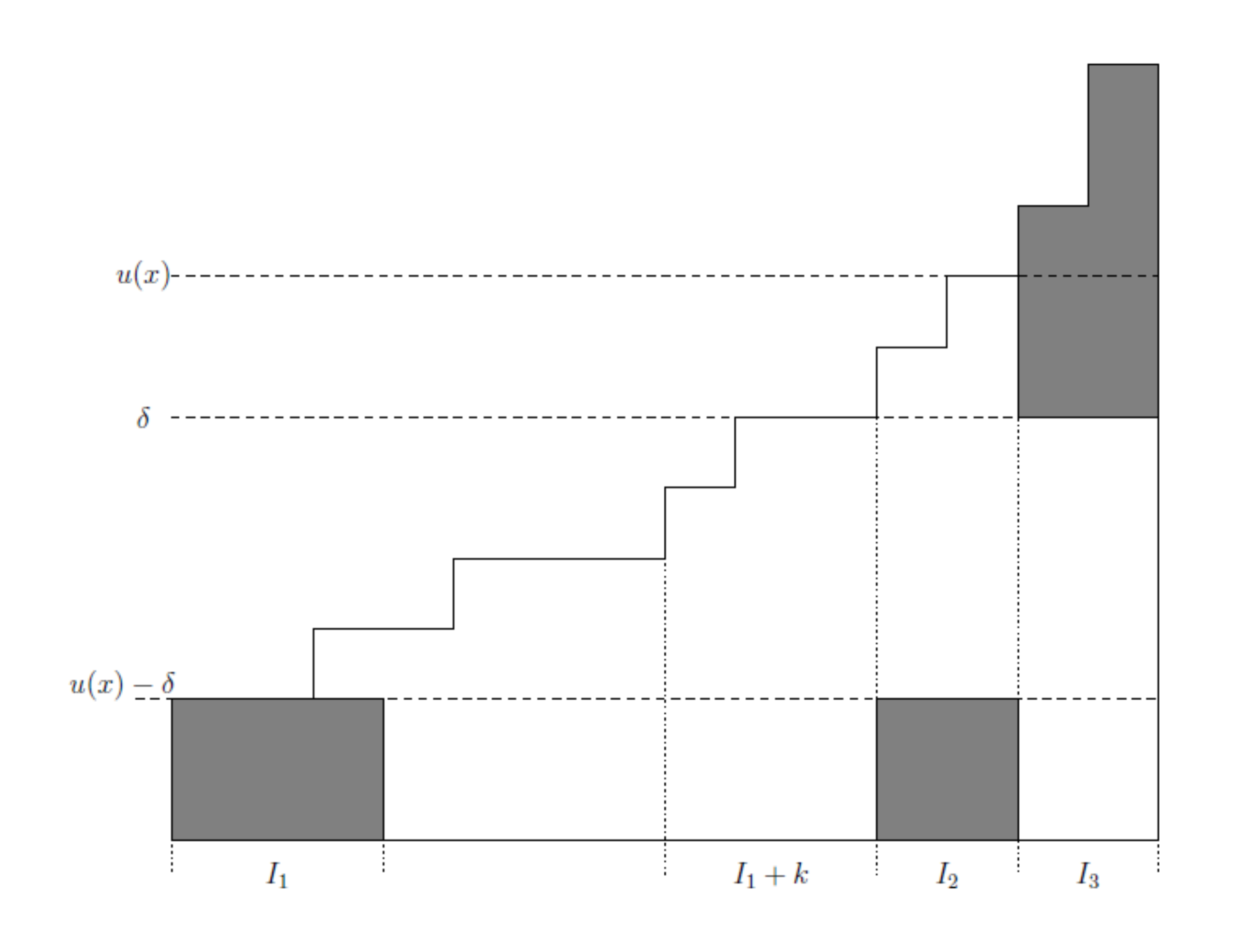}
\caption{$x'$ is obtained by removing the gray area.} \label{F3}
\end{figure}

\end{enumerate}

\subsubsection{Type II positions: $m(x)\geq z(x)$}

For a position $x$ let
\begin{align}
D_1(x) &=\{(m,y)\mid 0\leq m\leq m(x) , m(x)-m\leq y< y(x) \},\\
D_2(x) &=\{(m,y)\mid m(x)-v(x)+z(x)-1\leq m <m(x) , y=y(x) \},\\
D(x) &=D_1(x)\cup D_2(x).
\end{align}
Note that $D_2(x)=\emptyset$ whenever $(m(x)-z(x))\mod (y(x)+1)$ is zero.

\begin{lemma}\label{existsmoveinD}
For every $x$ and for every $(m,y) \in D(x)$ there is a position $x'$ reachable from $x$ such that $m(x')=m$ and $y(x')=y$.
\end{lemma}

\proof
For $(m,y) \in D(x)$ let us consider the following three cases.

Case 1:
$0 \leq m< m(x)$ and $x_2-m \leq y \leq y(x)$.
We consider two positions $\bar{x}$ and $\hat{x}$ reachable from $x$ defined as follows
$$\bar{x}_i =
\begin{cases}
m,    & \textrm{ if either } i=1 \textrm{ or } i\geq k+2;\\
x_i,  & \textrm{ otherwise};
\end{cases}
$$
$$
\hat{x}_i =
\begin{cases}
x_i,                    & \textrm{ if } i \leq k+1; \\
\min(x_i,m(x)+y(x))-1,  & \textrm{ if } i \geq k+2;
\end{cases}
$$
Since we have $y(\bar{x})=x_2-m$ and $y(\hat{x})\geq y(x)-1+m(x)-m \geq y(x)$,
we can apply Lemma \ref{continuity} and obtain that
for all values $x_2-m\leq y\leq y(x)$ there exists a position $x'$ such that
$\bar{x}\leq x'\leq \hat{x}$ and $y(x')=y$. All these positions have $m(x')=m$.

Let us also note that all these positions $x'$
are reachable from $x$, because exactly $k$ components of $x$ are decreased.

\medskip

Case 2: $0\leq m < m(x)$ and $m(x)-m \leq y < x_2-m$.

We consider two positions $\bar{x}$ and $\hat{x}$ reachable from $x$ defined as
$$
\bar{x}_i =
\begin{cases}
m,   & \textrm{ if } i=2 \textrm{ or } i\geq k+2;\\
x_i, & \textrm{ otherwise;}
\end{cases}
$$
$$\hat{x}_i =
\begin{cases}
m,   & \textrm{ if } i=2;\\
x_i, & \textrm{ if } i=1 \textrm{ or } i=3,4, \ldots, k+1; \\
\min(x_i,m(x)+y(x))-1, & \textrm{ if } i \geq k+2.
\end{cases}
$$

Since we have $y(\bar{x})=m(x)-m$ and
$y( \hat{x}) \geq x_{k+1}-m-1 \geq x_2-m-1$, we can apply Lemma \ref{continuity}, and obtain
that for all $y \in[m(x)-m,x_2-m-1]$ there exists a position $x' \in [\bar{x},\hat{x}]$ such that $y=y(x')$.

\medskip

Case 3: $m=m(x)$ and $0 \leq y < y(x)$.
Let us note that the last inequaliti implies $y(x) \geq 1$ and therefore $x_{k+1}>m(x)$.
We consider two positions $\bar{x}$ and $\hat{x}$ reachable from $x$ defined as
$$\bar{x}_i =
\begin{cases}
x_i,   & \textrm{ if } i=1, \ldots ,k;\\
m(x),  & \textrm{ if } i\geq k+1; \\
\end{cases}
~~~~\hat{x}_i =
\begin{cases}
x_i,   & \textrm{ if } i=1, \ldots ,k; \\
x_i-1, & \textrm{ if } i\geq k+1. \\
\end{cases}
$$
Similarly to the previous cases, we have $y(\bar{x})=0$ and
$y(\hat{x}) = y(x)-1$, and thus by Lemma \ref{continuity} it follows that
for all $y \in[0,y(x)-1]$ there exists an $x' \in [\bar{x},\hat{x}]$ such that $y=y(x')$.

If we put all three cases together we cover all values $(m,y) \in D$.
\qed

\medskip

Let us set
\begin{equation} \label{v(m,y)}
v(m,y)= {y+1 \choose 2}+ \Big[ \Big(m-1-{y+1 \choose 2}\Big)\mod (y+1)\Big].
\end{equation}
Note that if $m=m(x)$ and $y=y(x)$ then $v(m,y)=v(x)$.
Furthermore, we have
$$V(y):=\{v(m,y) | m \in \ZZ_\geq\}=\left[ {y+1 \choose 2},{y+2 \choose 2} \right).$$
Therefore the sets $V(y)$, $y \in \ZZ_\geq$ form a partition of $\ZZ_\geq$, as shown in Lemma \ref{interval}.

\begin{lemma}\label{Dstar}
If $m(x) \geq z(x)$, then every $(m,y) \in D(x)$ satisfies the following relations
$$m\geq {y+1 \choose 2}+1~~\textrm{ and }~~\{v(m,y) \mid (m,y) \in D(x) \} = [0,v(x)).$$
\end{lemma}

\proof
Let us first consider $(m,y) \in D_1(x)$.
By the definition of $D_1$ and the assumption of  $m(x) \geq z(x)={y(x)+1 \choose 2}+1$ we get
$$m\geq m(x)-y \geq {y(x)+1 \choose 2}+1-y \geq {y+2 \choose 2}+1-y \geq {y+1 \choose 2}+1.$$
By the definition of $D_1(x)$, for all $y \in [0,y(x))$ we have $(m,y) \in D_1(x)$ for all $m \in [m(x)-y,m(x)]$. Hence, by (\ref{v(m,y)}) we have $\{v(m,y) \mid m \in [m(x)-y,m(x)]\}=V(y)$.
Thus, by Lemma~\ref{interval} we get $$\{v(m,y) \mid (m,y) \in D_1(x) \}=\bigcup_{y \in [0,y(x))}V(y)=[0,z(x)-1).$$

\smallskip

Let us next consider $(m,y) \in D_2(x)$.
By the definition of $ v(x)$ we can write $v(x)=z(x)-1+r $, where $r=m(x)-z(x)-\lambda(y(x)+1)$ for some $\lambda \in \ZZ_\geq$. This implies $m(x)-r\geq z(x)$. By the definition of $D_2(x)$ we have $m\geq m(x)-r$.
These two inequalities imply $m\geq {y+1 \choose 2} +1$.
Since $m \in [m(x)-r, m(x))$ takes $r$ consecutive values, we have
$\{v(m,y) \mid (m,y) \in D_2(x) \}=[z(x)-1,v(x))$.
\qed

The above lemma implies that any position $x'$ with $m(x')=m$, $y(x')=y$ for some $(m,y)\in D$
is a type II position. Hence $g(x')=v(x')=v(m,y)$.
Thus the second claim in the above lemma together with Lemma \ref{existsmoveinD} implies that for any $0\leq \delta <v(x)$ there exists a move $x\to x'$ such that $g(x')=\delta$.

Since we proved this for both type I and type II positions we concluded the proof of (II).

Properties (I) and (II) together now imply that $\G=g$.
This concludes the proof of Theorem \ref{thm.n=2k}.
\qed

\smallskip


\begin{figure}[ht]
\centering
\includegraphics[height=7cm]{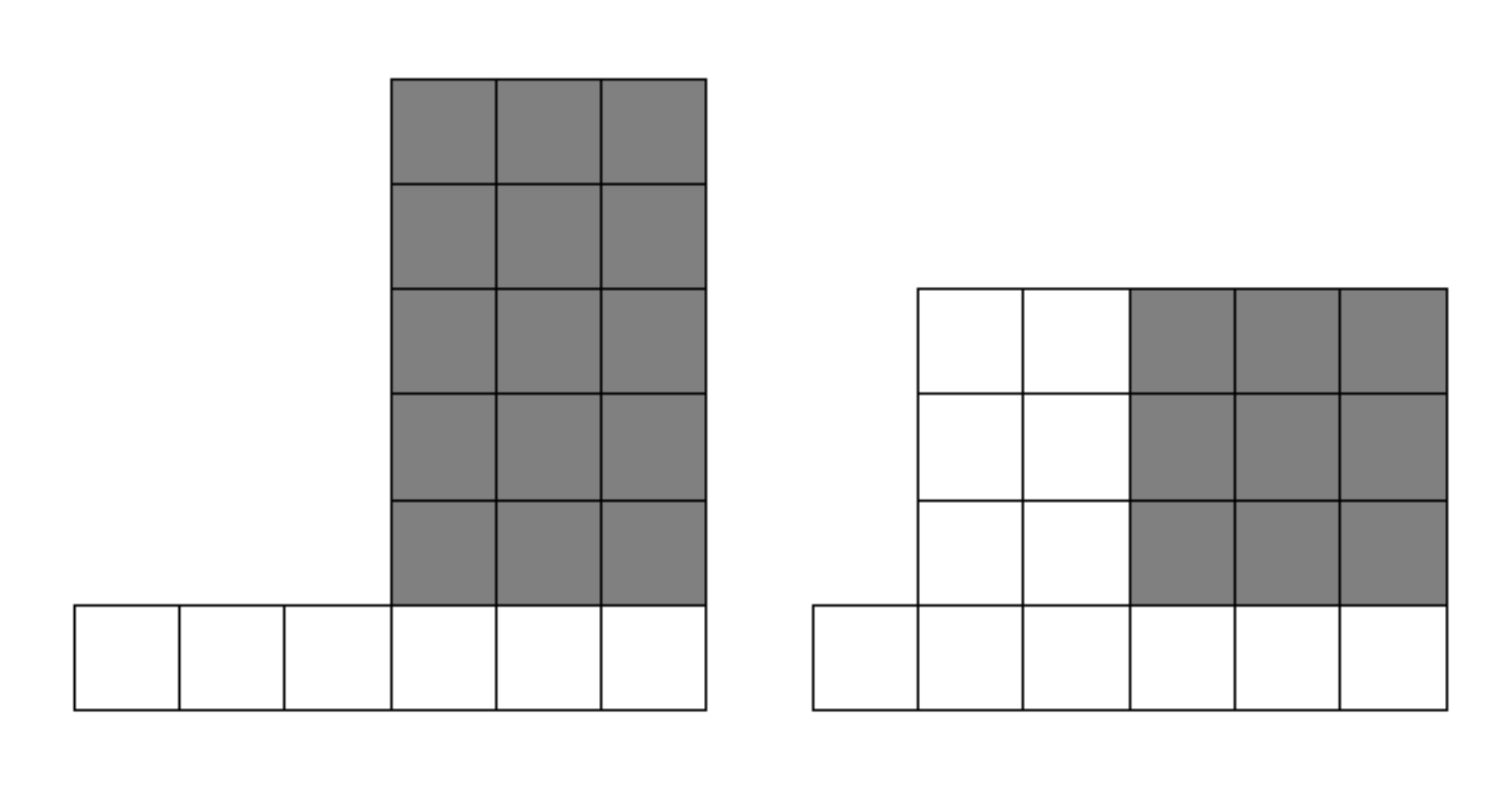}
\caption{Case $n=6,k=3$. Both above positions have $m=1,u=7,y=5$; but the position obtained by removing the gray area has $m=1,u=2,y=0$ and such a position cannot be obtained by the picture on the right since $m=1$ implies that $u \geq 4$.} \label{n=6,k=3}
\end{figure}

\section{Characterization of the $0$- and $1$-positions in Moore's game} \label{Ss.Moore01}

The following result was shown in \cite{JM80}. We provide here a different proof for the convenience of the reader.

\begin{theorem} $($see \cite{JM80}$)$. \label{t-Moore-m=0,1}
For any $k \in \{1, \ldots, n\}$  and  $m \in \{0,1\}$, a position
$x$  is an $m$-position of the Moore game  {\sc Nim}$^\leq_{n, k}$
if and only if  $M(x) = m$.
\end{theorem}

Let us represent the components of $x$ as binary sequences
\[
x_i=\sum_{j=0}^\infty x_{ij}2^j \text{ for } i=1,...,n,
\]
and define
\begin{equation}\label{e0001}
y_j=\sum_{i=1}^n x_{ij} \mod (k+1) \text{ for } j=0,1,....
\end{equation}
Then we have $M(x)=\sum_{j=0}^\infty y_j(k+1)^j$.

\begin{lemma}\label{moore1}
Let  $x\to x'$ be a move, and let $j$ be the highest index such that $x_{ij}\neq x'_{ij}$ for some $i$.
Then we must have $x'_{ij}\leq x_{ij}$ for all $i=1,...,n$.
\end{lemma}
\proof
In a move we can only decrease the components of $x$.
Therefore, if $x'_{ij}>x_{ij}$, then we must have a $j'>j$ such that $x'_{ij'}<x_{ij'}$.
\qed

\subsection{Proof of Theorem \ref{t-Moore-m=0,1} for $m=0$}   \label{ss-p-t1-0}

First, let us prove Moore's result:
$\cG(x)=0$ if and only if  $M(x)=0$.

\smallskip

By the properties of $\P$-positions it is enough to show that
\begin{itemize}
\item [\rm{(i0)}] for any position $x$ with $M(x)=0$, there exists no  move $x \to x'$ such that $M(x') = 0$;
\item [\rm{(a0)}] for any position $x$ with $M(x)> 0$,  there exists a move $x \to x'$ such that  $M(x') = 0$.
\end{itemize}

To show \rm{(i0)}, let us consider a move $x \to x'$
from a position  $x$  with  $M(x)=0$.

Let  $j$  be the highest binary bit such that
$x_{ij}$ and $x'_{ij}$  differ for some $i$. Such a $j$ must exist since in a move we must change at least one components.
By Lemma \ref{moore1} we have $x'_{ij}\leq x_{ij}$ for all $i=1,...,n$,
implying $1\leq \sum_{i=1}^n(x_{ij}-x'_{ij}) \leq k$ because  in a move we can change at most $k$ components. Therefore,
$(\sum_{i}x'_{ij} \mod (k+1)) \not= 0$ and, thus,  $M(x') > 0$.

\medskip

To show \rm{(a0)}, let us consider a position  $x$  with  $M(x)>0$. We will construct a move $x \to x'$ such that $M(x') = 0$.

\begin{notation} \label{T.N}
Let  $t_1, \dots , t_p$ denote the bits $j$ such that $y_j\not=0$, assuming  $t_1 > \dots > t_p$.
Set $N = \{1,2,\ldots,n\}$.
\end{notation}

The following algorithm defines index sets $\emptyset=I_0\subseteq I_1\subseteq \cdots I_p \subseteq N$
such that we can compute a move $x \to x'$ with $M(x')=0$, by decreasing components $i\in I_p$. Define $O_j=\{i\in I_j\mid x_{it_{j+1}}=1\}$, and set $\ga_j = |I_j|$  and  $\gb_j=|O_j|$.

\medskip
\noindent
{\bf Step 0}. Initialize $I_0=\emptyset$
and hence, $\ga_0 = \gb_0 = 0$, and set $x'_i:=x_i$ for all $i \in N$.

\smallskip
\noindent
{\bf Step 1}. For $j=1, \dots , p$, construct $I_j$ and update $x'$ as follows.

\smallskip

{\bf Case 1.} If $y_{t_j} \leq \beta_{j-1}$, then let  $I_j:=I_{j-1}$,
choose  $y_{t_j}$ many indices $i\in O_{j-1}$,
and update $x'_{it_j}:=0$.

\smallskip

{\bf Case 2.} If $y_{t_j} > \beta_{j-1}$ and $(k+1) -y_{t_j} \leq \ga_{j-1}-\gb_{j-1}$,  then
let  $I_j:=I_{j-1}$, choose  $(k+1)-y_{t_j}$ many indices $i\in I_{j-1}\setminus O_{j-1}$, and update $x'_{i{t_j}}:=1$.

\smallskip

{\bf Case 3.}  If $y_{t_j} > \beta_{j-1}$ and $y_{t_j}-\beta_{j-1} \leq k-\ga_{j-1}$, then
let $I_j$ be the index set obtained from $I_{j-1}$ by adding
$y_{t_j} - \beta_{j-1}$ many indices $i$ from $N \setminus I_{j-1}$ such that
$x_{i{t_j}}=1$.  Update $x'_{i{t_j}}:=0$ for $i\in (I_j\setminus I_{j-1})\cup O_{j-1}$.

\medskip
\noindent
We first note that the three cases above are exclusive and
cover all possible $y_{t_j}$, $\ga_{j-1}$ and $\gb_{j-1}$ values.
Moreover, it is easily seen that
a position $x'$ after the execution of the algorithm satisfies  $M(x')=0$, and
$x'_i=x_i$  holds for $i \not\in I_p$.

Note that we increase the set $I_j$ only in Case 3, in which case we have
\[
|I_j|=\alpha_j=(y_{t_j}-\beta_{j-1})+\alpha_{j-1} \leq k,
\]
implying $|I_p|\leq k$.

Note also that in Case 3 we must have at least $y_{t_j} - \beta_{j-1}$ many indices $i\in N\setminus I_{j-1}$ with $x_{it_j}=1$ by the definition of $y_{t_j}$ in \eqref{e0001}.

It remains to show that $x' < x$.
Assume that $i$ is an index such that $x'_i \not= x_i$.
Then some $j$ satisfies $i \not\in I_{j-1}$ and $i \in I_j$.
This implies that $x'_i$ was first updated during the $j$th iteration of Step 1.
Namely, the $t_j$th bit of $x'_i$ is modified from $1$ to $0$.
Since $x'_{it}=x_{it}$ holds for all $t$ with $t > t_j$,
we have $x'_i < x_i$, which completes the proof.
\qed

\subsection{Proof of Theorem \ref{t-Moore-m=0,1} for $m=1$}
\label{ss-p-t1-1}

Now, let us prove that
$\cG(x)=1$ if and only if  $M(x)=1$.

\smallskip
The proof in the previous subsection implies that for a position $x$ with $M(x)=1$ there exists a move $x\to x'$ such that $M(x')=0$.
By the properties of the SG function, it remains to show that
\begin{itemize}
\item [\rm{(i1)}] for any position $x$  with $M(x) = 1$, there exists
no move $x \to x'$ such that $M(x') = 1$;
\item [\rm{(a1)}] for any position $x$ with $M(x) > 1$, there
exists a move $x \to x'$ such that  $M(x')=1$;
\end{itemize}

We prove \rm{(i1)} similarly to (i0). 
Let us assume that $M(x)=1$ holds for a position $x$
and consider a move $x \to x'$.
Let $j$ be the highest binary bit such that
$x_{ij}$ and $x'_{ij}$ differ for some $i$.
Then by Lemma \ref{moore1} we have $x'_{ij} \leq x_{ij}$ for all $i$,  and
$1 \leq \sum_{i}(x_{ij}-x'_{ij}) \leq k$.
Hence, $\sum_{i}x'_{ij} \not= \sum_{i}x_{ij}  (\mod (k+1))$ and, thus, $M(x') \not= 1$.

To show \rm{(a1)}, let us consider a position $x$  with $M(x) > 1$.
Similarly to  \rm{(a0)}, we will algorithmically construct a move $x \to x'$ such that $M(x')=1$.

Let again $t_1, \dots , t_{p-1}$ denote the bits $j>0$ such that $y_j>0$,
where we assume that  $t_1 > \dots > t_{p-1}$, and add $t_p=0$.

The algorithm remains the same, as for \rm{(a0)},
except for $j=p$, when $t_p=0$. We detail below the computation of $I_p$ from $I_{p-1}$:

\medskip
\noindent

{\bf Case 1.} If $y_0> 1$ and $y_0-1 \leq \gb_{p-1}$,
then let  $I_p:=I_{p-1}$,
choose  $y_0-1$ many indices ${i}$ from $I_p$  such that
$x'_{i0}=1$, and update $x'_{i0}:=0$ for such indices.

\smallskip

{\bf Case 2.} If $y_0> 1$ and $(k+2)-y_0 \leq \ga_{p-1}-\beta_{p-1}$,
let  $I_p:=I_{p-1}$, choose  $(k+2)-y_{0}$ many indices ${i}$ from $I_p$  such
that $x'_{i0}=0$ and update $x'_{i0}:= 1$ for such indices.

\smallskip

{\bf Case 3.} If $y_0> 1$ and $0< (y_0-1) -\gb_{p-1} \leq k-\ga_{p-1}$, then
let $I_p$ be an index set obtained from $I_{p-1}$ by adding $(y_{0}-1)-\beta_{p-1}$ many
indices $i$ from $N \setminus I_{p-1}$ such that $x_{i{0}}=1$, and
update $x'_{i0}:=0$ for all $i \in I_p$ with $x'_{i0}=1$.

\smallskip

{\bf Case 4.} If $y_0=0$ and $\ga_{p-1} > \gb_{p-1}$,
then let  $I_p:=I_{p-1}$,
choose  an index ${i}$ from $I_p$  such that $x'_{i0}=0$, and update $x'_{i0}:=1$.

\smallskip

{\bf Case 5.} If $y_0=0$ and $\ga_{p-1} = \gb_{p-1}=k$,
then let  $I_p:=I_{p-1}$, and update $x'_{i0}:=0$ for all $i \in I_p$.

\smallskip

{\bf Case 6.} If $y_0=0$ and $\ga_{p-1} = \gb_{p-1}<k$,
then let $I_p$ be an index set obtained from $I_{p-1}$ by adding $k-\ga_{p-1}$ many
indices $i$ from $N \setminus I_{p-1}$ such
that $x_{i{0}}=1$, and update $x'_{i0}:=0$ for all $i \in I_p$.

\smallskip
{\bf Case 7.} If $y_0=1$, then we set $I_p=I_{p-1}$.

\medskip
\noindent
Note that the above seven cases are exclusive and cover all possible $y_{0}$, $\ga_{p-1}$, and $\gb_{p-1}$ values.
Note also that
in Case 6,  $\ga_{p-1} > 0$ since otherwise  $M(x)=0$, giving a contradiction.
Thus, in Case 6, we can choose $k - \ga_{p-1}$ many
indices $i$ from $N \setminus I_{p-1}$ such that $x_{i{0}}=1$.
Let $x'$ be a position obtained by the algorithm.
Then, clearly  $M(x')=1$  and  $x \to x'$ is a move. This completes the proof of (a1).
\qed

\section{More on the Tetris function} \label{Ss.Tetris}

In this section we show that the SG function described in Theorems \ref{thm.n<2k}
and \ref{thm.n=2k} can in fact be computed efficiently, in polynomial time. For this we need to show that the Tetris function can be computed in polynomial time for these games.
We also prove that for a given position $x$ and integer $0\leq g < T_{n,k}(x)$
we can compute in polynomial time a move $x\to x'$ such that $T_{n,k}(x')=g$.
Finally, in subsection \ref{tetris-degree} we recall
some relations to degree sequences of graphs and hypergraphs.

\subsection{Computing the Tetris function in polynomial time} \label{comp-tet-in-p}

Let us recall that to a position $x$ we associated a shifted position $\bar x$ after Corollary \ref{delta-epsilon} with the property that $T_{n,k}(x)=\bar x_{n-k+1}$.
The procedure described there is a non-polynomial algorithm.
However $\bar x$ and consequently $T_{n,k}(x)=\bar x_{n-k+1}$ can be computed in a more efficient way.

\begin{theorem}\label{Tetrispolynomial}
Given a position $x$ we can compute $T_{n,k}(x)$ in linear time in $n$.
\end{theorem}

\proof
We can assume without loss of generality that $x$ is a nondecreasing position. We show that the corresponding $\bar x$ can be constructed in linear time in $n$, and thus the claim follows by the equality $T_{n,k}(x)=\bar x_{n-k+1}$.

Recall that the input size is $\log(\prod_{i = 1}^n x_i)$. Let $s=\sum_{i=1}^{n-k} x_i$ be the number of tokens we shift on top of the largest $k$ piles; see Figure \ref{ADD(1)}.
We know that for some $\ell <k$ the first $\ell +1$ columns of $\bar x$ have almost the same number of tokens (at most one difference.)
To determine this index $\ell$ and the height of the resulting piles, we use simple volume based arguments.
We need to compute first the following parameters.

For each $i=1, \ldots ,k-1$, we denote by $y_i=x_{n-k+i+1}-x_{n-k+i}$ the difference of the sizes of consecutive piles.
Set $s_0 =0$, $s_{k}=\infty$, and
for $i=1,...,k-1$, set $s_i=s_{i-1}+i \cdot y_i$ (i.e., the number of tokens we need to shift on top of the first $i$ piles $(n-k+1), \ldots, (n-k+i)$ to make them all equal to $x_{n-k+1+i}$.)
We define a unique $\ell$ by $s_\ell \leq s < s_{\ell+1}$.
We define $a = s - s_\ell$, $\alpha = \lfloor\frac{a}{\ell+1}\rfloor$, and $\beta =a \mod (\ell+1)$.
We fill up the first $\ell+1$ columns to level $x_{n-k+\ell+1}$ using $s_\ell$ tokens. Then,
we place the remaining $a$ tokens by increasing each of the first
$\ell+1$ columns (indexed $n-k+1,...,n-k+\ell +1$) by $\alpha$ and the last $\beta$ of these by one more, as in the following expression.
\begin{figure}[ht]
\centering
\includegraphics[height=8cm]{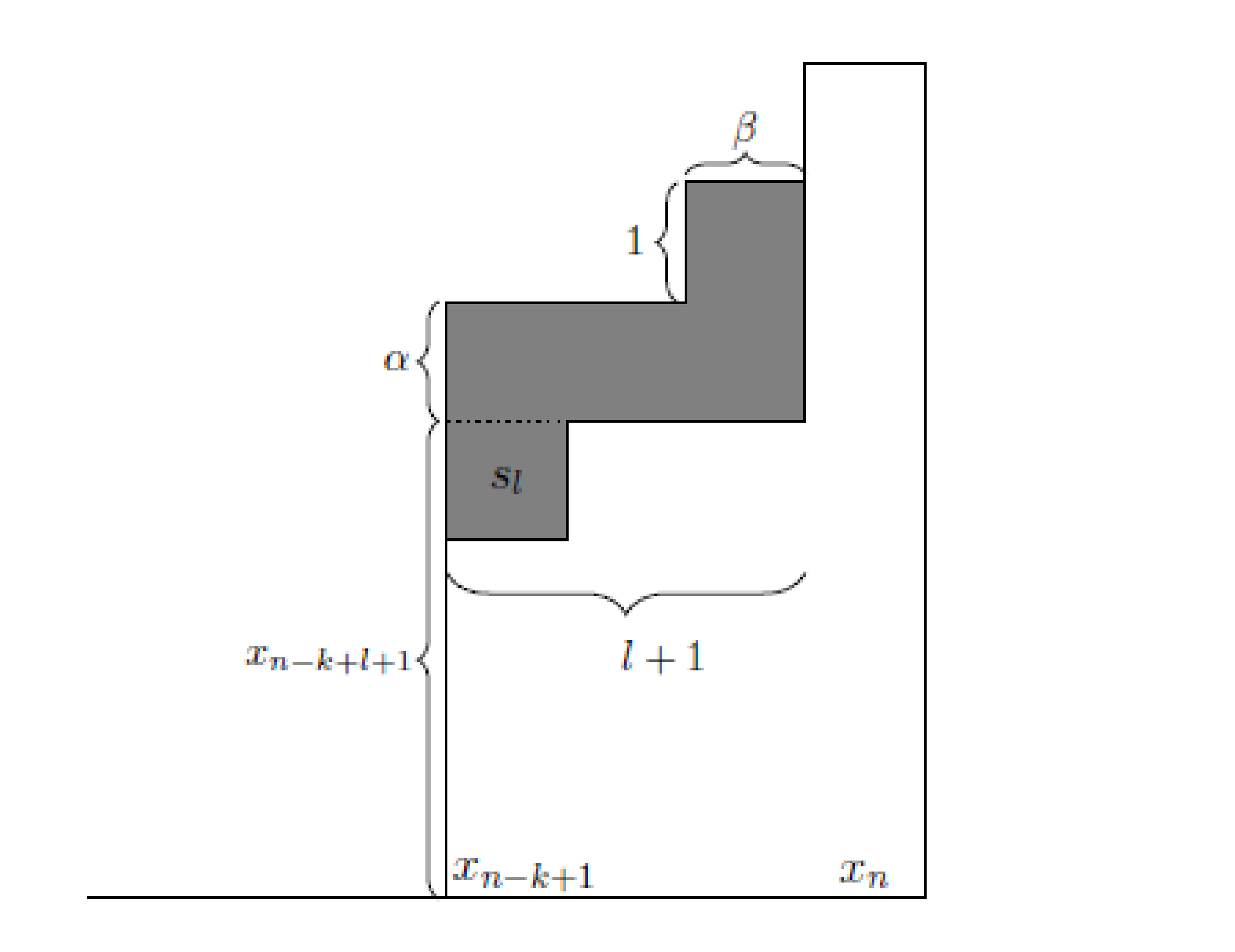}
\caption{An example of calculating $\bar x$ for $x=(1,2,2,3,4,4,7)$ with $k = 4$.}
\end{figure}

\begin{align*}
\bar x_i=
\begin{cases}
0,                        & \text{if } i=1, \ldots ,n-k; \\
x_{n-k+\ell+1}+\alpha,    & \text{if } i=n-k+1, \ldots ,n-k+1+\ell -\beta; \\
x_{n-k+\ell+1}+\alpha +1, & \text{if } i=n-k+2+\ell -\beta, \ldots, n-k+1+\ell; \\
x_i,                      & \text{if } i=n-k+2+\ell, \ldots ,n.
\end{cases}
\end{align*}
It is easy to see that this defines $\bar x$ correctly,
and that all these parameters can be computed in $O(n)$ time, if $x$ is a nondecreasing vector.
\qed

\begin{remark}
Technically, only the computation of $s=\sum_{i=1}^{n-k}x_i$ depends on $n$.
All other computations in the previous proof can be done in $O(k)$ time.
\end{remark}
\medskip

\subsection{Polynomial computation of a move to a given Tetris value} \label{}

Let $x = (x_1, \ldots, x_n)$ be a nondecreasing position and $k > n/2$.
We examine the question of how to move to some position of given Tetris value $g$.

We denote by
\begin{enumerate} \itemsep0em
\item [\rm{(i)}] $x^\ell$ the position obtained from $x$ by removing all tokens from the largest $k$ piles of $x$, and by
\item [\rm{(ii)}] $x^u$ the position obtained from $x$ by decreasing the largest $k$ piles by one unit each.
\end{enumerate}

Consider the set
$W = \{T_{n,k}(x') \mid \forall x': x\to x'\}$.
By Lemma \ref{tbasic} (i) we have $T_{n,k}(x^\ell) = \min(W)$.
By Lemma \ref{bestslowmove} we have $T_{n,k}(x^u)=T_{n,k}(x)-1$. As in the proof of
Theorem \ref{thm.n<2k}, we can argue that for
every value $g$ such that $T_{n,k}(x^\ell)\leq g \leq T_{n,k}(x)-1$ there exists a move $x\to x'$ such that  $x^\ell \leq x' \leq x$ and $T_{n,k}(x')=g$.

\begin{theorem}\label{achieveTetrisvalue}   
Given $g \in W$, computing a position $x'$ such that $x^\ell \leq x' \leq x$ and $T_{n,k}(x')=g$ can be done in $O\big(n\log(\sum_{i=n-k+1}^nx_i)\big)$ time.
\end{theorem}
\proof
We have $T_{n,k}(x^\ell)\leq g \leq T_{n,k}(x^u)$.
Using the monotonicity of the Tetris function we perform a binary search in the space of positions between $x^\ell$ and $x^u$.
In a general step we compute $L=\sum_{i=n-k+1}^nx^\ell_i$ and $U=\sum_{i=n-k+1}^nx^u_i$, set $M=\lfloor\frac{L+U}{2}\rfloor$ and compute $y_i=int(\frac{x^\ell_i+x^u_i}{2})$ for $i=n-k+1$,
where $int(\cdot)$ is a rounding to a nearest integer value in such a way that $\sum_{i=n-k+1}^ny_i=M$.
Finally, we set $y_i=x_i$ for $i<n-k+1$.
If $T_{n,k}(y)<g$ then we replace $x^\ell$ by $y$, otherwise we replace $x^u$ by $y$.

Clearly these computations can be done in each step in $O(n)$ time,
and computing the Tetris value of $y$ can also be done in $O(n)$ time by Theorem~\ref{Tetrispolynomial}.
\qed

\smallskip

\begin{remark}
Similarly to the proof of Proposition \ref{Tetrispolynomial}
we could improve the complexity of the above algorithm to $O(n)$.
\end{remark}

\subsection{Tetris function and degree sequences of graphs and hypergraphs}\label{tetris-degree}

A related problem is the hypergraph realization of a given degree sequence.
Let us fix $V=\{1,2,...,n\}$ as the set of vertices.
A multi-hypergraph $\cH=\{H_1,...,H_m\}$ is a family of subsets (called hyperedges) of $V$,
i.e., $H_j\subseteq V$ for all $j=1,...,m$. It is called $k$-uniform if $|H_j|=k$ for all $j=1,...,m$.
The degree $d_{\cH}(i)$ of a vertex $i\in V$ is the number of hyperedges $H_j$ of $\cH$ that contain $i$.
We allow the same subset to appear multiple times in $\cH$.

Given an integer vector $x\in\ZZ_{\geq}^n$, one can ask if there exists a $k$-uniform multi-hypergraph $\cH$ on the vertex set $V$ such that $d_{\cH}(i)=x_i$ for all $i\in V$.
Equivalently, we examine the existence of a bipartite graph $G=(X,Y,E)$ such that $|X|=n$, $|Y|=m$, $d_G(i)=x_i$ for all $i\in X$, and $d_G(j)=k$ for all $j\in Y$.
For the latter we can apply the classical Gale-Ryser theorem claiming that the answer is yes if and only if

\begin{itemize} \itemsep0em
\item [\rm{(i)}]  $\sum_{i=1}^nx_i=km$
\item [\rm{(ii)}] $\sum_{i=1}^n\min(x_i,g) \geq kg$ for all $g=1,...,m$.
\end{itemize}
Let us note that checking these conditions may not be polynomial in $x$ and $k$, since $m=\sum_{i=1}^nx_i/k$ according to property \rm{(i)}.
Let us also note that following a sequence of slow moves starting from position $x$, each time the set of columns that we decrease by $1$ can be considered as a hyperedge of $\cH$.
Thus a maximal sequence of slow moves will construct $\cH$ if the Tetris function achieves its trivial bound $kT_{n,k}(x)=\sum_{i=1}^nx_i$.
This equality is in fact equivalent with property \rm{(i)}, since we must have $m=T_{n,k}(x) $ in this case.
Our results in this section thus prove that for the above degree sequence realization problems the most efficient answer is to compute the Tetris function value in linear time, and then compare it to its trivial upper bound.
If these are the same then the answer is yes.

Havel (1955) and Hakimi (1962) provided a simple greedy algorithm based on a characterization for the recognition of degree sequences of bipartite graphs.
For the above case their criterion states that $x$ is a degree sequence of a $k$-uniform multi-hypergraph if and only if the position $x'$ is also a degree sequence of a $k$-uniform multi-hypergraph, where $x'$ is obtained from $x$ by decreasing the $k$ largest components of $x$ by $1$.
Note that this implies a recursive process that is one of the definitions we used for the Tetris function.

Let us remark finally that in general $x$ is not the degree sequence of a $k$-uniform multi-hypergraph.
In this case however a move $x\to x'$ such that $T_{n,k}(x')=0$ provides us with a minimal modification such that $x''=x-x'$ becomes the degree sequence of such a hypergraph.

\end{document}